\documentclass{amsart}
\usepackage{amssymb,latexsym}

\theoremstyle{plain}
\newtheorem{theorem}{Theorem}
\newtheorem{corollary}[theorem]{Corollary}

\newtheorem{proposition}[theorem]{Proposition}

\theoremstyle{definition}
\newtheorem{definition}[theorem]{Definition}

\newtheorem*{remark}{Remark}

\def \cL {{\mathcal L}}
\def \cF {{\mathcal F}}

\def \sdu {\bigtriangleup}
\def \sdd {\bigtriangledown}

\def \ff {
\unitlength .4 mm
\linethickness{0.4pt}
\begin{picture}(12,10)(0,2.5)
\put(0.00,0.00){\circle*{0.5}}
\put(10.00,0.00){\circle*{0.5}}
\put(0.00,10.00){\circle*{0.5}}
\put(10.00,10.00){\circle*{0.5}}
\put(5.00,5.00){\circle*{0.5}}
\put(5.00,8.0){\circle*{2.5}}
\put(2.0,5.00){\circle*{2.5}}
\put(8.0,5.00){\circle*{2.5}}
\put(5.00,2.0){\circle*{2.5}}
\put(0,0){\line(1,0){10}}
\put(0,0){\line(0,1){10}}
\put(10,0){\line(0,1){10}}
\put(0,10){\line(1,0){10}}
\end{picture}
}

\def \ee {
\unitlength .4 mm
\linethickness{0.4pt}
\begin{picture}(12,10)(0,2.5)
\put(0.00,0.00){\circle*{0.5}}
\put(10.00,0.00){\circle*{0.5}}
\put(0.00,10.00){\circle*{0.5}}
\put(10.00,10.00){\circle*{0.5}}
\put(5.00,5.00){\circle*{0.5}}
\put(5.00,8.0){\circle{2.5}}
\put(2.0,5.00){\circle{2.5}}
\put(8.0,5.00){\circle{2.5}}
\put(5.00,2.0){\circle{2.5}}
\end{picture}
}

\def \x {
\unitlength .4 mm
\linethickness{0.4pt}
\begin{picture}(12,10)(0,2.5)
\put(0.00,0.00){\circle*{0.5}}
\put(10.00,0.00){\circle*{0.5}}
\put(0.00,10.00){\circle*{0.5}}
\put(10.00,10.00){\circle*{0.5}}
\put(5.00,5.00){\circle*{0.5}}
\put(5.00,8.0){\circle{2.5}}
\put(2.0,5.00){\circle*{2.5}}
\put(8.0,5.00){\circle{2.5}}
\put(5.00,2.0){\circle*{2.5}}
\put(0,0){\line(1,0){10}}
\put(0,0){\line(0,1){10}}
\end{picture}
}

\def \y {
\unitlength .4 mm
\linethickness{0.4pt}
\begin{picture}(12,10)(0,2.5)
\put(0.00,0.00){\circle*{0.5}}
\put(10.00,0.00){\circle*{0.5}}
\put(0.00,10.00){\circle*{0.5}}
\put(10.00,10.00){\circle*{0.5}}
\put(5.00,5.00){\circle*{0.5}}
\put(5.00,8.0){\circle{2.5}}
\put(2.0,5.00){\circle{2.5}}
\put(8.0,5.00){\circle*{2.5}}
\put(5.00,2.0){\circle*{2.5}}
\put(0,0){\line(1,0){10}}
\put(10,0){\line(0,1){10}}
\end{picture}
}

\def \xs {
\unitlength .4 mm
\linethickness{0.4pt}
\begin{picture}(10.00,10.00)(0,2.5)
\put(0.00,0.00){\circle*{0.5}}
\put(10.00,0.00){\circle*{0.5}}
\put(0.00,10.00){\circle*{0.5}}
\put(10.00,10.00){\circle*{0.5}}
\put(5.00,5.00){\circle*{0.5}}
\put(5.00,8.0){\circle*{2.5}}
\put(2.0,5.00){\circle{2.5}}
\put(8.0,5.00){\circle*{2.5}}
\put(5.00,2.0){\circle{2.5}}
\put(10,0){\line(0,1){10}}
\put(0,10){\line(1,0){10}}
\end{picture}
}

\def \xdy {
\unitlength .4 mm
\linethickness{0.4pt}
\begin{picture}(12,10)(0,2.5)
\put(0.00,0.00){\circle*{0.5}}
\put(10.00,0.00){\circle*{0.5}}
\put(0.00,10.00){\circle*{0.5}}
\put(10.00,10.00){\circle*{0.5}}
\put(5.00,5.00){\circle*{0.5}}
\put(5.00,8.0){\circle{2.5}}
\put(2.0,5.00){\circle*{2.5}}
\put(8.0,5.00){\circle*{2.5}}
\put(5.00,2.0){\circle{2.5}}
\put(0,0){\line(1,0){10}}
\put(0,0){\line(0,1){10}}
\put(10,0){\line(0,1){10}}
\put(0,10){\line(1,0){10}}
\end{picture}
}

\def \xny {
\unitlength .4 mm
\linethickness{0.4pt}
\begin{picture}(12,10)(0,2.5)
\put(0.00,0.00){\circle*{0.5}}
\put(10.00,0.00){\circle*{0.5}}
\put(0.00,10.00){\circle*{0.5}}
\put(10.00,10.00){\circle*{0.5}}
\put(5.00,5.00){\circle*{0.5}}
\put(5.00,8.0){\circle{2.5}}
\put(2.0,5.00){\circle*{2.5}}
\put(8.0,5.00){\circle*{2.5}}
\put(5.00,2.0){\circle{2.5}}
\end{picture}
}

\def \xply {
\unitlength .4 mm
\linethickness{0.4pt}
\begin{picture}(12,10)(0,2.5)
\put(0.00,0.00){\circle*{0.5}}
\put(10.00,0.00){\circle*{0.5}}
\put(0.00,10.00){\circle*{0.5}}
\put(10.00,10.00){\circle*{0.5}}
\put(5.00,5.00){\circle*{0.5}}
\put(5.00,8.0){\circle{2.5}}
\put(2.0,5.00){\circle*{2.5}}
\put(8.0,5.00){\circle*{2.5}}
\put(5.00,2.0){\circle{2.5}}
\put(0,0){\line(1,0){10}}
\put(0,0){\line(0,1){10}}
\end{picture}
}

\def \xpry {
\unitlength .4 mm
\linethickness{0.4pt}
\begin{picture}(12,10)(0,2.5)
\put(0.00,0.00){\circle*{0.5}}
\put(10.00,0.00){\circle*{0.5}}
\put(0.00,10.00){\circle*{0.5}}
\put(10.00,10.00){\circle*{0.5}}
\put(5.00,5.00){\circle*{0.5}}
\put(5.00,8.0){\circle{2.5}}
\put(2.0,5.00){\circle*{2.5}}
\put(8.0,5.00){\circle*{2.5}}
\put(5.00,2.0){\circle{2.5}}
\put(0,0){\line(1,0){10}}
\put(10,0){\line(0,1){10}}
\end{picture}
}

\def \xlsy {
\unitlength .4 mm
\linethickness{0.4pt}
\begin{picture}(12,10)(0,2.5)
\put(0.00,0.00){\circle*{0.5}}
\put(10.00,0.00){\circle*{0.5}}
\put(0.00,10.00){\circle*{0.5}}
\put(10.00,10.00){\circle*{0.5}}
\put(5.00,5.00){\circle*{0.5}}
\put(5.00,8.0){\circle{2.5}}
\put(2.0,5.00){\circle*{2.5}}
\put(8.0,5.00){\circle*{2.5}}
\put(5.00,2.0){\circle{2.5}}
\put(10,0){\line(0,1){10}}
\put(0,10){\line(1,0){10}}
\end{picture}
}

\def \xrsy {
\unitlength .4 mm
\linethickness{0.4pt}
\begin{picture}(12,10)(0,2.5)
\put(0.00,0.00){\circle*{0.5}}
\put(10.00,0.00){\circle*{0.5}}
\put(0.00,10.00){\circle*{0.5}}
\put(10.00,10.00){\circle*{0.5}}
\put(5.00,5.00){\circle*{0.5}}
\put(5.00,8.0){\circle{2.5}}
\put(2.0,5.00){\circle*{2.5}}
\put(8.0,5.00){\circle*{2.5}}
\put(5.00,2.0){\circle{2.5}}
\put(0,0){\line(0,1){10}}
\put(0,10){\line(1,0){10}}
\end{picture}
}

\def \xat {
\unitlength .4 mm
\linethickness{0.4pt}
\begin{picture}(12,10)(0,2.5)
\put(0.00,0.00){\circle*{0.5}}
\put(10.00,0.00){\circle*{0.5}}
\put(0.00,10.00){\circle*{0.5}}
\put(10.00,10.00){\circle*{0.5}}
\put(5.00,5.00){\circle*{0.5}}
\put(5.00,8.0){\circle{2.5}}
\put(2.0,5.00){\circle{2.5}}
\put(8.0,5.00){\circle{2.5}}
\put(5.00,2.0){\circle{2.5}}
\put(0,0){\line(1,0){10}}
\put(0,0){\line(0,1){10}}
\end{picture}
}

\def \xsat {
\unitlength .4 mm
\linethickness{0.4pt}
\begin{picture}(12,10)(0,2.5)
\put(0.00,0.00){\circle*{0.5}}
\put(10.00,0.00){\circle*{0.5}}
\put(0.00,10.00){\circle*{0.5}}
\put(10.00,10.00){\circle*{0.5}}
\put(5.00,5.00){\circle*{0.5}}
\put(5.00,8.0){\circle{2.5}}
\put(2.0,5.00){\circle{2.5}}
\put(8.0,5.00){\circle{2.5}}
\put(5.00,2.0){\circle{2.5}}
\put(10,0){\line(0,1){10}}
\put(0,10){\line(1,0){10}}
\end{picture}
}

\def \xb {
\unitlength .4 mm
\linethickness{0.4pt}
\begin{picture}(10.00,10.00)(0,2.5)
\put(0.00,0.00){\circle*{0.5}}
\put(10.00,0.00){\circle*{0.5}}
\put(0.00,10.00){\circle*{0.5}}
\put(10.00,10.00){\circle*{0.5}}
\put(5.00,5.00){\circle*{0.5}}
\put(5.00,8.0){\circle{2.5}}
\put(2.0,5.00){\circle*{2.5}}
\put(8.0,5.00){\circle{2.5}}
\put(5.00,2.0){\circle*{2.5}}
\end{picture}
}

\def \xsvy {
\unitlength .4 mm
\linethickness{0.4pt}
\begin{picture}(10.00,10.00)(0,2.5)
\put(0.00,0.00){\circle*{0.5}}
\put(10.00,0.00){\circle*{0.5}}
\put(0.00,10.00){\circle*{0.5}}
\put(10.00,10.00){\circle*{0.5}}
\put(5.00,5.00){\circle*{0.5}}
\put(5.00,8.0){\circle*{2.5}}
\put(2.0,5.00){\circle{2.5}}
\put(8.0,5.00){\circle*{2.5}}
\put(5.00,2.0){\circle*{2.5}}
\put(0,0){\line(1,0){10}}
\put(0,0){\line(0,1){10}}
\put(10,0){\line(0,1){10}}
\put(0,10){\line(1,0){10}}
\end{picture}
}

\def \xsgy {
\unitlength .4 mm
\linethickness{0.4pt}
\begin{picture}(10.00,10.00)(0,2.5)
\put(0.00,0.00){\circle*{0.5}}
\put(10.00,0.00){\circle*{0.5}}
\put(0.00,10.00){\circle*{0.5}}
\put(10.00,10.00){\circle*{0.5}}
\put(5.00,5.00){\circle*{0.5}}
\put(5.00,8.0){\circle{2.5}}
\put(2.0,5.00){\circle{2.5}}
\put(8.0,5.00){\circle*{2.5}}
\put(5.00,2.0){\circle{2.5}}
\end{picture}
}

\def \xgys {
\unitlength .4 mm
\linethickness{0.4pt}
\begin{picture}(10.00,10.00)(0,2.5)
\put(0.00,0.00){\circle*{0.5}}
\put(10.00,0.00){\circle*{0.5}}
\put(0.00,10.00){\circle*{0.5}}
\put(10.00,10.00){\circle*{0.5}}
\put(5.00,5.00){\circle*{0.5}}
\put(5.00,8.0){\circle{2.5}}
\put(2.0,5.00){\circle*{2.5}}
\put(8.0,5.00){\circle{2.5}}
\put(5.00,2.0){\circle{2.5}}
\end{picture}
}

\def \xsvys {
\unitlength .4 mm
\linethickness{0.4pt}
\begin{picture}(10.00,10.00)(0,2.5)
\put(0.00,0.00){\circle*{0.5}}
\put(10.00,0.00){\circle*{0.5}}
\put(0.00,10.00){\circle*{0.5}}
\put(10.00,10.00){\circle*{0.5}}
\put(5.00,5.00){\circle*{0.5}}
\put(5.00,8.0){\circle*{2.5}}
\put(2.0,5.00){\circle*{2.5}}
\put(8.0,5.00){\circle*{2.5}}
\put(5.00,2.0){\circle{2.5}}
\put(0,0){\line(1,0){10}}
\put(0,0){\line(0,1){10}}
\put(10,0){\line(0,1){10}}
\put(0,10){\line(1,0){10}}
\end{picture}
}

\def \bool {
\unitlength .4 mm
\linethickness{0.4pt}
\begin{picture}(10.00,10.00)(0,2.5)
\put(0.00,0.00){\circle*{0.5}}
\put(10.00,0.00){\circle*{0.5}}
\put(0.00,10.00){\circle*{0.5}}
\put(10.00,10.00){\circle*{0.5}}
\put(5.00,5.00){\circle*{0.5}}
\put(5.00,8.0){\circle*{2.5}}
\put(2.0,5.00){\circle*{2.5}}
\put(8.0,5.00){\circle*{2.5}}
\put(5.00,2.0){\circle*{2.5}}
\end{picture}
}

\def \xbysat {
\unitlength .4 mm
\linethickness{0.4pt}
\begin{picture}(10.00,10.00)(0,2.5)
\put(0.00,0.00){\circle*{0.5}}
\put(10.00,0.00){\circle*{0.5}}
\put(0.00,10.00){\circle*{0.5}}
\put(10.00,10.00){\circle*{0.5}}
\put(5.00,5.00){\circle*{0.5}}
\put(5.00,8.0){\circle{2.5}}
\put(2.0,5.00){\circle*{2.5}}
\put(8.0,5.00){\circle{2.5}}
\put(5.00,2.0){\circle*{2.5}}
\put(0,0){\line(0,1){10}}
\put(0,10){\line(1,0){10}}
\end{picture}
}

\def \xatyb {
\unitlength .4 mm
\linethickness{0.4pt}
\begin{picture}(10.00,10.00)(0,2.5)
\put(0.00,0.00){\circle*{0.5}}
\put(10.00,0.00){\circle*{0.5}}
\put(0.00,10.00){\circle*{0.5}}
\put(10.00,10.00){\circle*{0.5}}
\put(5.00,5.00){\circle*{0.5}}
\put(5.00,8.0){\circle{2.5}}
\put(2.0,5.00){\circle{2.5}}
\put(8.0,5.00){\circle*{2.5}}
\put(5.00,2.0){\circle*{2.5}}
\put(0,0){\line(1,0){10}}
\put(0,0){\line(0,1){10}}
\end{picture}
}

\def \xatxyb {
\unitlength .4 mm
\linethickness{0.4pt}
\begin{picture}(10.00,10.00)(0,2.5)
\put(0.00,0.00){\circle*{0.5}}
\put(10.00,0.00){\circle*{0.5}}
\put(0.00,10.00){\circle*{0.5}}
\put(10.00,10.00){\circle*{0.5}}
\put(5.00,5.00){\circle*{0.5}}
\put(5.00,8.0){\circle{2.5}}
\put(2.0,5.00){\circle{2.5}}
\put(8.0,5.00){\circle{2.5}}
\put(5.00,2.0){\circle*{2.5}}
\put(0,0){\line(1,0){10}}
\put(0,0){\line(0,1){10}}
\end{picture}
}

\begin{document}

\title{Non--Commutative Symmetric Differences in Orthomodular Lattices}

\author{Gerhard Dorfer}

\address{ Institute of Algebra and Computational Mathematics\\ 
Vienna University of Technology\\
Wiedner Hauptstr.\ 8-10/118, A-1040 Vienna, Austria}
\email{g.dorfer@tuwien.ac.at}

\keywords{Orthomodular lattice, symmetric difference, congruence relation}
\subjclass{Primary: 06\,C\,15; Secondary: 06\,B\,10, 81\,P\,10}
\thanks{Research supported by \"OAD, Cooperation between Austria and Czech Republic in Science and Technology, 
grant no.\ 2003/1, and by the Austrian Science Fund FWF under project no.\ S~8312.} 

\begin{abstract}
We deal with the following question: What is the proper way to introduce 
symmetric difference in orthomodular lattices? 
Imposing two natural conditions on this operation, six possibilities remain: 
the two (commutative) normal forms of the symmetric difference
in Boolean algebras and four non-commutative terms. 
It turns out that in many respects the non-commutative forms, though more complex with respect to the
lattice operations, in their properties are much nearer to the symmetric difference in Boolean algebras
than the commutative terms. As application we demonstrate the usefulness of non-commutative
symmetric differences in the context of congruence relations.
\end{abstract}

\maketitle

\markboth{G.~Dorfer}{Non-Commutative Symmetric Differences in OML}

\section{Introduction of symmetric differences}
The symmetric difference plays a prominent role in the theory of Boolean algebras (BA).
For instance, important properties of congruence relations in BA such as permutability, 
regularity and uniformity of congruences follow mainly from the fact that the symmetric difference
is an associative, cancellative and invertible term function. We will recall all these notions later in detail when
we deal with it.

Thus it is a manifest task to investigate symmetric difference in the more 
general framework of orthomodular lattices (OML). Some work in this direction can be found in \cite{DDL}
and \cite{La}.

An orthomodular lattice $\cL=(L,\lor,\land,',0,1)$ is a bounded 
lattice $(L,\lor,\land,0,1)$ with an orthocomplementation $'$, i.e., for all 
$x,y\in L$ 
\[
  x\land x'=0,\;x\lor x'=1,\;x''=x,\;x\leq y\
  \mbox{ implies }\ y'\leq x',
\]
and $\cL$ satisfies the orthomodular law:
\[ x\leq y\ \mbox{ implies }\ y=x\lor (y\land x').\] 

As distinguished from Boolean algebras orthomodular lattices are not distributive. 
The following two relations provide some kind of a measure for non-dis\-tri\-bu\-tivity in a particular OML:
\begin{itemize}
  \item[--] The {\em commutativity relation} $C$: $aCb$ if and only if the subalgebra generated by
  $\{ a,b\}$ in $\cL$ is Boolean. For instance, $a\leq b$ or $a\leq b'$ imply
  $aCb$.
  \item[--] The {\em perspectivity relation} $\sim$: $a\sim b$ if and only if $a$ and $b$ have a
  common (algebraic) complement, i.e.\ there exists an element $c\in L$ such
  that  
$$a\land c=b\land c=0,\quad a\lor c=b\lor c=1.$$
\end{itemize}
An OML $\cL$ is a BA if and only if $C$ is the all~relation or,
equivalently, if and only if $\sim$ is the identity. 

For a solid introduction to the theory of OML we refer to \cite{Ka}.

We recall a smart technique of Navara \cite{Na} to represent elements and simplify computations in the free
OML $\cF(x,y)$ with free generators $x$ and $y$: Let $c(x,y):=(x\land y)\lor (x\land y')\lor
(x'\land y)\lor(x'\land y')$ denote the commutator of $x$ and $y$. Instead of $(c(x,y))'$ we simply write $c'(x,y)$.
$\cF(x,y)\cong [0,c(x,y)]\times [0,c'(x,y)]$, where $[0,c(x,y)]\cong 2^4$ is
the 16-element BA (which is the free BA generated by two elements) with atoms $x\land y, x\land y', x'\land y,
x'\land y'$, and $[0,c'(x,y)]\cong\mathrm{MO2}$ is the six-element OML with atoms 
$x\land c'(x,y),y\land c'(x,y),x'\land c'(x,y),y'\land c'(x,y)$. 

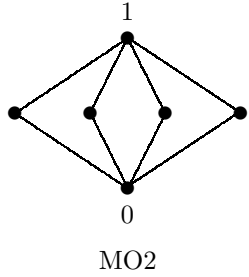
\unitlength 0.5mm
\linethickness{0.4pt}
\begin{picture}(60.00,80.00)(-52,75)
\multiput(70.00,100.00)(-0.18,0.12){167}{\line(-1,0){0.18}}
\multiput(40.00,120.00)(0.18,0.12){167}{\line(1,0){0.18}}
\multiput(70.00,140.00)(0.18,-0.12){167}{\line(1,0){0.18}}
\multiput(99.33,120.00)(-0.18,-0.12){167}{\line(-1,0){0.18}}
\multiput(70.00,100.00)(-0.12,0.24){84}{\line(0,1){0.24}}
\multiput(60.00,120.00)(0.12,0.24){84}{\line(0,1){0.24}}
\multiput(70.00,140.00)(0.12,-0.24){84}{\line(0,-1){0.24}}
\multiput(80.00,120.00)(-0.12,-0.24){84}{\line(0,-1){0.24}}
\put(70.00,100.00){\line(0,1){0.00}}
\put(70.00,100.00){\circle*{3.33}}
\put(40.00,120.00){\circle*{3.33}}
\put(60.00,120.00){\circle*{3.33}}
\put(80.00,120.00){\circle*{3.33}}
\put(100.00,120.00){\circle*{3.33}}
\put(70.00,140.00){\circle*{3.33}}
\put(70.00,147.00){\makebox(0,0)[cc]{$1$}}
\put(70.00,93.00){\makebox(0,0)[cc]{$0$}}
\put(70.00,81.00){\makebox(0,0)[cc]{MO2}}
\end{picture}

\noindent 
The representation of elements in 
$\cF(x,y)$ refers to the following scheme:

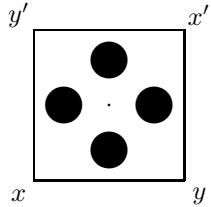
\begin{center}
\unitlength .2 cm
\linethickness{0.5pt}
\begin{picture}(12.00,13.00)(-1,-1)
%\put(0.00,0.00){\circle*{0.5}}
%\put(10.00,0.00){\circle*{0.5}}
%\put(0.00,10.00){\circle*{0.5}}
%\put(10.00,10.00){\circle*{0.5}}
\put(5.00,5.00){\circle*{0.05}}
\put(5.00,8.0){\circle*{2.5}}
\put(2.0,5.00){\circle*{2.5}}
\put(8.0,5.00){\circle*{2.5}}
\put(5.00,2.0){\circle*{2.5}}
\put(0,0){\line(1,0){10}}
\put(0,0){\line(0,1){10}}
\put(10,0){\line(0,1){10}}
\put(0,10){\line(1,0){10}}
\put(-1,-1){\makebox(0,0)[cc]{$x$}}
\put(11,-1){\makebox(0,0)[cc]{$y$}}
\put(-1,11){\makebox(0,0)[cc]{$y'$}}
\put(11,11){\makebox(0,0)[cc]{$x'$}}
\end{picture}
\end{center}

\noindent
The discs correspond to the Boolean part and the bars to the MO2-part of $\cF(x,y)$.
Full/empty discs refer to the presence/absence of the corresponding atoms in the Boolean part and
corner angles represent the atoms in the MO2-part. 
While we deal with $\cF(x,y)$, deviating from \cite{Na}, we do not indicate the generators
$x$ and $y$ separately. In the representation one just has to remember that $x$ is the element
down left, $y$ down right and complements are vis-\`{a}-vis. For instance, 
$$(x\land y')\lor(x'\land y)=\xny,\quad x'\land c'(x,y)=\xsat,\quad 1=\ff.$$
Computations in $\cF(x,y)$ decompose into a Boolean part with set-theoretical operations on the discs, and
a MO2-part with operations on the corner angles following the evaluation rules in MO2. For instance,
$$
\begin{array}{l}
x=(x\land c'(x,y))\lor((x\land y)\lor(x\land y'))=\xat\lor\xb=\x\, ,\\ 
\ \\
x'\lor y=\left(\x\right)'\lor \y=\xs\lor\y =\xsvy=\left( \xgys\right)'=\left(x\land y'\right)'.
\end{array}$$
In an attempt to adopt symmetric difference for OML the first striking thing is that two
different terms representing the symmetric difference in a BA may differ when they are evaluated
in an OML. In particular, consider the disjunctive and conjunctive normal form
$$
\begin{array}{c}
    x\sdd y:=(x\land y')\lor (x'\land y),\\
    x\sdu y:=(x\lor y)\land (x'\lor y').
  \end{array}
$$
Applying these operations in MO2 with generating elements $a,b$,
we obtain
$$
  a\sdd b=a\sdd b'=a'\sdd b=a'\sdd b'=0,\quad a\sdu b=a\sdu b'=a'\sdu b=
  a'\sdu b'=1.
$$
In fact the difference between these two operations could not be larger.

First of all we have to make clear what we understand by a symmetric difference. We impose the following 
two natural conditions.
\begin{definition}\label{defsd}
A binary operation $+$ in an OML $\cL$ is called  {\em symmetric difference} if
\begin{itemize}
\item[(i)] $+$ is a {\em term function},
\item[(ii)] $+$ coincides with the conventional symmetric difference if $\cL$ is a BA.
\end{itemize}
\end{definition}
%\noindent 
The first objective is to find out how many such operations exist.
\begin{theorem}\label{thm1}
For OML there are exactly six possibilities to define an operation such that
(i) and (ii) are satisfied:
$$
\begin{array}{rcl}
x \bigtriangledown y & = & (x \land y')\lor (x' \land y) \\
x\bigtriangleup y & = & (x\lor y)\land (x'\lor y') \\
x +_l y & := & (x \lor (x' \land y))\land (x'\lor y') \\
x+_r y & := & ((x\land y')\lor y)\land (x'\lor y') \\
x+_{l'}y & := & (x\lor y)\land(x'\lor(x\land y')) \\
x+_{r'}y & := & (x\lor y)\land((x'\land y)\lor y')
\end{array}
$$
\end{theorem}
\begin{proof}
Let $\cF_{\rm BA}(u,v)$ denote the free BA with free generators
$u$ and $v$.
We consider the homomorphism
$$
    \varphi:\left\{
    \begin{array}{ccc}
      \cF(x,y) & \to & \cF_{\rm BA} (u,v)\\
      x & \mapsto & u\\
      y & \mapsto & v
    \end{array}\right. .
$$
Condition (i) and (ii) in Definition~\ref{defsd} exactly mean that the symmetric differences 
are given by the terms in $\varphi^{-1}(u\sdu v)$. 
Using Navara's technique these elements have the following representation:
$$\begin{array}{l}
\xny =(x\land y')\lor(x'\land y)=x\sdd y\\
\xdy =((x\land y')\lor(x'\land y))\lor c'(x,y)=(x\lor y)\land (x'\lor y')=x\sdu y\\
\xply=((x\land y')\lor(x'\land y))\lor(x\land c'(x,y))=(x\lor(x'\land y))\land(x'\lor y')\\
\xpry=((x\land y')\lor(x'\land y))\lor(y\land c'(x,y))=((x\land y')\lor y)\land(x'\lor y')\\
\xlsy=((x\land y')\lor(x'\land y))\lor(x'\land c'(x,y))=(x\lor y)\land(x'\lor(x\land y'))\\
\xrsy=((x\land y')\lor(x'\land y))\lor(y'\land c'(x,y))=(x\lor y)\land((x'\land y)\lor y')
\end{array}
$$
\end{proof}
%\noindent
To point up the difference between the six terms we rewrite two of them in Navara's notation:
$$ x+_l y=\x +_l \y =\xply,\qquad x\sdu y=\x\sdu \y=\xdy.$$
In the Boolean part in both cases the conventional symmetric difference is formed, in the MO2-part
the arguments do not commute and $+_l$ results in the {\em left} argument whereas $\sdu$ results in
the join.

Next we summarize some properties of symmetric differences.
We omit the proof which is straightforward.
\begin{proposition}\label{p1}
In $\cF(x,y)$ the following holds true.
\begin{enumerate} 
\item[\rm (i)] $x\sdu y=y\sdu x=x'\sdu y'=(x\lor y)\sdu(x\land y)$,
\item[\rm (ii)] $x\sdd y=(x\sdu y')'$, $x\sdu y=(x\sdd y')'$,
\item[\rm (iii)] $x\sdu y$ and $x\sdd y$ commute with both $x$ and $y$,
\item[\rm (iv)] $x+_l y=y+_r x=x'+_{l'}y'=y'+_{r'}x'$,
\item[\rm (v)] $(x+_l y)'=x'+_l y$,
\item[\rm (vi)] $x+_l y$ commutes with $x$ but does not commute with $y$.
\end{enumerate}
\end{proposition}
%\noindent
In this proposition it becomes evident that there are strong interrelations between the six symmetric differences 
and that they naturally split into two subclasses: one consisting of the commutative 
$\sdu$ and $\sdd$ and the other one containing the remaining four non-commutative terms.
Operations within the same subclass behave very similar and, as will turn out, operations from
different classes differ in their properties. Thus in the following we will state results
for $\sdu$ and $+_l$ only. These results may be reformulated for the other symmetric differences by
the help of Proposition~\ref{p1} easily.

Let in the following $\cL=(L,\lor,\land,',0,1)$ denote an arbitrary OML 
and $a,b,c$ elements of $L$. 
The next proposition will clarify to what extent the six symmetric differences may differ.
\begin{proposition}
\begin{enumerate}
\item[\rm (i)] If $a$ commutes with $b$ then all 
six symmetric differences of $a$ and $b$ are equal.
\item [\rm (ii)] If $a$ does not commute with $b$ 
then the six symmetric differences of $a$ and $b$
are pairwise different.
\end{enumerate}
\end{proposition}
\begin{proof}
In $\cF(x,y)$ the six symmetric differences of $x$ and $y$ form an interval isomorphic to MO2,
in $\cL$ the symmetric differences of $a$ and $b$ thus form a homomorphic image of MO2. 
Since MO2 is simple, this image either
consists of one element, which is the case if $a$ commutes with $b$, 
or is isomorphic to MO2 if $a$ does not commute with $b$.
\end{proof} 
\begin{corollary}
\begin{enumerate}
\item[\rm (i)]
If two distinct symmetric differences coincide on the whole of $\cL$ then $\cL$ is a BA.
\item[\rm (ii)]
If $+_l$ is commutative then $\cL$ is a BA.
\end{enumerate}
\end{corollary}
\section{Cancellativity, invertibility and associativity of symmetric differences}
In this section we study symmetric differences in OML with respect to 
important properties the symmetric difference fulfils in a BA.

For the convenience of the reader we recall some basic notions from algebra.
A binary operation $\circ$ on a set $A$ is called right cancellative (left invertible)
%see E.S. Ljapin, Semigroups, Translations of Mathemtical Monographs, Vol. 3, AMS, Providence, Rhode Island, 1963
%cancellativity widely used (Google search)
if for arbitrary $a,b\in A$ the equation $x\circ a=b$ has at most (at least) one solution
$x\in A$. Left cancellativity (right invertibility) is defined accordingly.
An operation is cancellative (invertible) if it is both left and right cancellative (left and right invertible).

Firstly we recall some known results (cf.\ \cite[Proposition~3.4, Lemma~3.6]{DDL}).
\begin{proposition}\label{p2}
\begin{enumerate}
\item[\rm (i)] Two elements $a$ and $b$ are perspective to each other
if and only if there exists $c$ such that $a\sdu c=b\sdu c.$
\item[\rm(ii)]  Two elements $a$ and $b$ commute if and only if there exists an element $c$
such that $a\sdu c=b$.
\end{enumerate}
\end{proposition}
%\noindent
By (i) we may interpret the relation of perspectivity in an OML as a measure for how
far the symmetric difference $\sdu$ is from being cancellative, and by (ii) the same
applies to the complement of the commutativity relation with respect to invertibility of $\sdu$.
\begin{theorem}\label{thm2}
In the variety of OML the symmetric difference $+_l$ satisfies the identity
$$ (x+_l y)+_l y=x.$$
\end{theorem}
\begin{corollary}
\begin{enumerate}
\item[\rm (i)] The symmetric difference $+_l$ is right cancellative, i.e., 
if $a+_l b=c+_l b$ then $a=c.$
\item[\rm (ii)] $+_l$ is left invertible, i.e.,
for all $a$ and $b$ there exists $c$ such that $c+_l a=b.$
%\item[\rm (iii)] $(L,+_l)$ is a right-quasigroup.
\end{enumerate}
\end{corollary}
\begin{proof}
Again we apply Navara's technique to verify the identity:
$$(x+_l y)+_l y=\left( \x +_l \y\right)+_l \y=\xply +_l \y=\x =x.$$
The corollary now follows easily. (i): If $a+_l b=c+_l b$ then by $+_l$-adding $b$ from the right 
we obtain $a=c.$

(ii): Choose $c=b+_l a$.
\end{proof}
In comparison to Proposition~\ref{p2} we characterize elementwise 
under which conditions left cancellation is
possible and a right inverse exists for the operation $+_l$.
\begin{proposition}
\begin{enumerate}
\item[\rm (i)] In general $+_l$ is not left cancellative, in particular $b+_l a=b+_l c$ 
if and only if $b\land a=b\land c$  and $b'\land a=b'\land c.$
\item[\rm (ii)] Two elements $a$ and $b$ commute if and only if
there exists $c$ such that $a+_l c=b.$
\end{enumerate}
\end{proposition}
%\noindent
Here we employ Navara's technique (as it was done in \cite{Na}) to represent terms in two variables
without assuming that these variables are free generators. As a consequence the representations may not
be unique in this case. To distinguish this from the case with free generators we indicate the generators as
indices in the corresponding positions. For example in MO2 with generators $a$ and $b$ we have
$c(a,b)=0$, i.e., ${}_{a\,}\bool_{\;b}={}_{a\,}\ee_b$.
\begin{proof}
(i): If $b+_l a=b+_l c$, i.e., $_{b\,}\xply_{a}= {}_{b\,}\xply_{c}$ then intersecting both sides with
$b'= {}_{b\,}\xs_{\,a}= {}_{b\,}\xs_{\,c}$ we arrive at $ {}_{b\,}\xsgy_{\,a}= {}_{b\,}\xsgy_{\,c}$,
i.e., $b'\land a=b'\land c$.

Joining both sides with $b'$ we obtain ${}_{b\,}\xsvys_{\,a}= {}_{b\,}\xsvys_{\,c}$,
hence $b'\lor a'=b'\lor c'$, i.e., $b\land a=b\land c$.

Conversely, if $b\land a=b\land c$ and $b'\land a=b'\land c$ then immediately 
$b+_l a=b+_l c$ follows.

(ii) is obvious.
\end{proof}
\begin{corollary}
If a symmetric difference is cancellative or invertible on $\cL$ then $\cL$ is a BA. 
\end{corollary}
%\noindent
In this context we recall \cite[Theorem~3.8]{DDL}.
\begin{proposition}
In the variety of OML there does not exist a binary term inducing a cancellative,
respectively invertible, term function on every OML.
\end{proposition}
As next step we address associativity of symmetric differences in OML.
\begin{proposition}
The following are equivalent:
\begin{enumerate}
\item[\rm (i)] $a$ commutes with $b$,
\item[\rm (ii)] $(a\bigtriangleup b)\bigtriangleup b=a\bigtriangleup(b\bigtriangleup b)$,
\item[\rm (iii)] $a+_l(a+_l b)=(a+_l a)+_l b$.
\end{enumerate}
\end{proposition}
\begin{proof}
Evidently (i) implies (ii) and (iii).

(ii) implies (i): Simplifying the left hand side we get $$(a\sdu b)\sdu b= {}_{a\,}\xdy_{\,b}\,\sdu\,
{}_{a\,}\y_{\,b}={}_{a\,}\xbysat_{\,b}.$$ The right hand side is $a\sdu(b\sdu b)=a=
{}_{a\,}\x_{\,b}$. Obviously ${}_{a\,}\xbysat_{\,b}={}_{a\,}\x_{\,b}$ implies that $a$ commutes with $b$.

A similar argument yields (iii) implies (i).
\end{proof}
%\noindent
We see that the associative law is valid for symmetric differences in OML in some special situations only.
The following corollary is derived easily.
\begin{corollary}
If a symmetric difference is associative in $\cL$ then $\cL$ is a BA. 
\end{corollary}
%\noindent
We also provide a positive result where again the discrepancy among symmetric differences appears.
\begin{proposition}
If $b$ commutes with $a$ and $c$ then
\begin{enumerate}
\item[\rm (i)]
$(a\sdu b)\sdu c=a\sdu(b\sdu c)$, 
\item[\rm (ii)]
$(b+_l a)+_l c=b+_l(a+_l c).$
\end{enumerate}
\end{proposition}
\begin{proof}
We proof (ii) in detail, (i) follows similarly.

Since $b$ commutes with $a$ and $c$ all computations occur in the commutator 
$C(b)\cong [0,b]\times [0,b']$ (see e.g.\ \cite[1.3.1.\ Theorem]{Ka}). Therefore it is sufficient to verify
the assertion for $b=0$ and $b=1$.

For $b=0$ we have 
$$(0+_l a)+_l c=a+_l c=0+_l(a+_lc),$$
and for $b=1$ due to Proposition~\ref{p1}
$$(1+_l a)+_l c=a'+_l c=(a+_l c)'=1+_l(a+_l c).$$
\end{proof}
%\noindent
Without going into detail we mention that also the distributivity of the meet operation
with respect to the symmetric difference cannot be generalized from BA to OML
even if one additionally considers also all {\em possible}\/ meet operations 
(in the sense of Definition~\ref{defsd}).
\section{Congruence relations}
It is  well-known that there is
a bijection between congruence relations of an OML $\cL$ and certain ideals of $\cL$,
so-called p-ideals \cite{Fi} (or orthomodular ideals \cite{Ka}): 
A lattice ideal $I$ is a p-ideal if it is closed under
perspectivity, i.e., $a\in I$ and $b\sim a$ imply $b\in I$. 
For a congruence $\theta$ on $\cL$ and $a\in L$ let $[a]\theta$ denote the congruence 
class of $a$. In the following theorem the relationship
between congruences and p-ideals is summarized.
\begin{theorem}
\begin{enumerate}
\item[\rm (i)] $I=[0]\theta$ for some congruence $\theta$ if and only if
$I$ is a p-ideal.
\item[\rm (ii)] A lattice ideal $I$ is a p-ideal if and only if for all $x$ in $L$
$$ x\land(I\lor x')\subseteq I$$
(where $x\land(I\lor x')=\{ x\land(i\lor x')\ |\ i\in I\}$).
\item[\rm (iii)] For a p-ideal $I$ the congruence $\theta$ corresponding to $I$ is given by the condition
$x\theta y$ if and only if $x\sdu y\in I.$
\end{enumerate}
\end{theorem}
%\noindent
Using a non-commutative symmetric difference this result can be modified
as follows.
\begin{proposition}
\begin{enumerate}
\item[\rm (ii')] A lattice ideal $I$ is a p-ideal if and only if for all $x$ in $L$
$$x+_l(I+_l x)\subseteq I.$$
\item[\rm (iii')] For a p-ideal $I$ the congruence $\theta$ corresponding to $I$ is given by the condition
$x\theta y$ if and only if $x+_l y\in I.$
\end{enumerate}
\end{proposition}
\begin{proof}
(ii'): We compare the elements $x\land(i\lor x')$ and $x+_l(i+_l x)$ for $x$ in $L$ and $i$ in $I$.
$$
\begin{array}{l}
x\land(i\lor x')={}_{x\,}\x_{\,i}\land\left( {}_{x\,}\y_{\,i}\lor {}_{x\,}\xs_{\,i}\right)=
{}_{x\,}\x_{\,i}\land {}_{x\,}\xsvy_{\,i}={}_{x\,}\xatxyb_{\,i},\\
x+_l(i+_l x)= {}_{x\,}\x_{\,i}+_l\left( {}_{x\,}\y_{\,i}+_l {}_{x\,}\x_{\,i}\right)=
{}_{x\,}\x_{\,i} +_l {}_{x\,}\xpry_{\,i}= {}_{x\,}\xatyb_{\,i}.  
\end{array}
$$
This means that $x+_l(i+_l x)=(x\land(i\lor x'))\lor(i\land x')$. Since $i\land x'$ is in $I$ 
($I$ is an order ideal) this means that the conditions in (ii) and (ii') are equivalent.

(iii'): We have to show that for a p-ideal $I$ the condition $x\sdu y\in I$ is equivalent to $x+_l y\in I$.
Since $x+_l y\le x\sdu y$ one direction is clear. On the other hand, if $x+_l y\in I$ then
$y+_l x=y+_l((x+_l y)+_l y)$ is also in $I$ (see (ii')), and $x\sdu y=(x+_l y)\lor(y+_l x)$. 
\end{proof}
\begin{remark} 
We want to point out that using non-commutative symmetric differences 
we can state the requirements for a
congruence kernel (p-ideal) of an OML very similar to groups/Boolean algebras.
One can even push this similarity further:
\begin{proposition}
A subset $I$ of $L$ is a congruence kernel  if and only if
\begin{enumerate}
\item[(i)] $(I,+_l,0)$ is a subalgebra of $(L,+_l,0)$ ({\em subgroup} condition),
\item[(ii)] $x+_l(I+_l x)\subseteq I$ for all $x\in L$ ({\em normal subgroup} condition),
\item[(iii)] $I$ is an order ideal in $(L,\le)$ ({\em as in BA}).
\end{enumerate}
\end{proposition}
%\noindent
Though these conditions seem very natural one has to be careful: for instance,
the second condition may {\em not}\/ be substituted by
$(x+_l I)+_l x\subseteq I$ for all $x\in L$. We leave the proof to the interested reader.
\end{remark}
Now we turn towards congruence classes. We recall some notions from universal algebra.
An algebra is called 
\begin{itemize}
\item[--] {\em congruence regular} if, for any congruence of this algebra, 
every congruence class determines the whole congruence uniquely,
\item[--] {\em congruence uniform} if for any fixed congruence all congruence classes have the same
cardinality,
\item[--] {\em congruence permutable} if any two congruences $\theta$ and $\phi$ permute, i.e.,
$$\theta\circ\phi=\phi\circ\theta\ (=\theta\lor\phi).$$   
\end{itemize}
A variety of algebras is called congruence regular (uniform, permutable) if all algebras of this variety
have this property.
For varieties congruence permutability and congruence regularity can be characterized by so-called
Mal'cev conditions (cf.\ \cite{Ma}, \cite{Cs}).
\begin{theorem}\label{thmmc}
\begin{enumerate}
\item[\rm (i)] 
A variety is congruence permutable if and only if there exists a ternary term $m(x,y,z)$ such that
the identities 
$$m(x,z,z)=x,\quad m(x,x,z)=z$$
hold true in the variety. $m(x,y,z)$ is called {\em Mal'cev term} for this variety.
\item[\rm (ii)] 
A variety is congruence regular if and only if there exist ternary terms 
$t_1(x,y,z),\ldots,t_n(x,y,z)$ such that
the condition
$$[ t_1(x,y,z)=z,\ldots,t_n(x,y,z)=z]\quad\mbox{if and only if}\quad x=y$$
is fulfilled in the variety. $t_1(x,y,z),\ldots,t_n(x,y,z)$ are called {\em Cs\'{a}k\'{a}ny terms} for this variety.
\end{enumerate}
\end{theorem}
\begin{theorem}\label{thm3}
In the variety of OML the term
$$(x+_l y)+_l z$$
is a Mal'cev and a Cs\'{a}k\'{a}ny term.
\end{theorem}
\begin{proof}
$(x+_l y)+_l z$ is a Mal'cev term:  By Theorem~\ref{thm2} we have
$$(x+_l z)+_l z=x,$$
and obviously
$$(x+_l x)+_l z=0+_l z=z.$$

$(x+_l y)+_l z$ is a Cs\'{a}k\'{a}ny term: The {\em if}\/ part of the condition in Theorem~\ref{thmmc}~(ii)
was verified just before. Now, suppose $(x+_l y)+_l z=z$, then by $+_l$-adding $z$ from the right side we get
$x+_l y=0$ and by adding $y$ we arrive at $x=y$.
\end{proof}
Let $I$ be a p-ideal, $\theta$ the congruence relation corresponding to $I$
and $a$ in $L$.
In \cite[Proposition~3.2,~3.3]{Do} the following relationship between 
the congruence classes of $\theta$ was given:
$$ [a]\theta=a\lor(I\land a'),\quad I=[a]\theta\sdu [a]\theta.$$
Using non-commutative symmetric differences much simpler (group-like) formulas
occur.
\begin{theorem}\label{thm4}
If $\theta$ is a congruence in $\cL$, $I=[0]\theta$ and $a\in L$ then
\begin{enumerate}
\item[\rm (i)] $[a]\theta=I+_l a,$
\item[\rm (ii)] $I=[a]\theta+_l a.$
\end{enumerate}
In particular, these formulas provide bijections between the congruence classes.
\end{theorem}
\begin{proof}
(i): If $i\in I$ then $(i+_l a)\theta(0+_l a)=a$. Conversely, if $x\in [a]\theta$ then $(x+_la)\theta
(a+_l a)=0$, hence $x+_l a\in I$ and by Theorem~\ref{thm2} $x=(x+_l a)+_l a\in I+_l a$.

(ii): $[a]\theta+_l a=(I+_l a)+_l a=I$. The concluding assertion follows easily.
\end{proof}
We should mention that with commutative symmetric differences no such 
simple formulas can be found. A similar phenomenon appears when 
terms representing implication in OML are investigated (cf.\ \cite{CHL}).

Finally we repeat a result (cf.\ e.g.\ \cite[Theorem~4.1--4.3]{Do}), which now is an immediate
consequence of Theorem~\ref{thm3} and Theorem~\ref{thm4}.
\begin{corollary}
The variety of OML is congruence regular, uniform and permutable.
\end{corollary}

\end{document}